\documentclass[a4paper,12pt,leqno]{article}
\usepackage{latexsym}
\usepackage[all]{xy}

\usepackage{amssymb} 
\usepackage{amsmath} 

\usepackage{tikz}
\usetikzlibrary{arrows,cd}  
\usepackage{amsthm}

\def\Z{{\mathbb{Z}}}
\def\K{{\mathbb{K}}}
\def\R{{\mathbb{R}}}
\def\C{{\mathbb{C}}}
\def\A{{\mathcal{A}}}
\def\B{{\mathcal{B}}}
\def\C{{\mathcal{C}}}

\DeclareMathOperator{\Der}{Der}

\DeclareMathOperator{\pd}{pd}

\DeclareMathOperator{\POexp}{POexp}

\numberwithin{equation}{section}

\newcommand{\owari}{\hfill$\square$}

\newtheorem{theorem}{Theorem}[section]
\newtheorem{prop}[theorem]{Proposition}
\newtheorem{cor}[theorem]{Corollary}
\newtheorem{lemma}[theorem]{Lemma}
\newtheorem{define}[theorem]{Definition}

\newtheorem{problem}[theorem]{Problem}

\theoremstyle{remark}

\newtheorem{example}[theorem]{Example}

\title{Generalization of the addition and restriction theorems from free arrangements to 
the class of projective dimension one}
\author{Takuro Abe}
\date{\today}

\begin{document}

\maketitle

\begin{abstract}
We study a generalized version of Terao's famous addition theorem for free arrangements to 
the category of those with projective dimension one. 
Namely, we give a criterion to determine the algebraic structure of logarithmic derivation 
modules of the addition when the deletion and restrictions are free with a mild condition. 
Also, we introduce a class of divisionally SPOG arrangements whose SPOGness depends only on the intersection lattice like Terao's famous conjecture on combinatoriality of freeness.
  \end{abstract}

\section{Introduction}
Let $\K$ be a field, $V=\K^\ell$, $S=\mbox{Sym}^*(V^*) =\K[x_1,\ldots,x_\ell]$ be the coordinate ring of 
$V$. For the derivation module $\Der S:=\oplus_{i=1}^\ell S\partial_{x_i}$ and a hyperplane arrangement 
$\A=\{H_1,\ldots,H_n\}$, where $H_i$ is defined as a zero locus of a non-zero linear form $\alpha_{H_i} \in V^*$, the \textbf{logarithmic derivation module} $D(\A)$ is defined by 
$$
D(\A):=\{\theta \in \Der S \mid \theta(\alpha_H) \in S \alpha_H\ (\forall H \in \A)\}.
$$
$D(\A)$ is an $S$-graded reflexive module of rank $\ell$, but not free in general. So we say that $\A$ is \textbf{free} with \textbf{exponents} $\exp(\A)=(d_1,\ldots,d_\ell)$ if $D(\A) \simeq \oplus_{i=1}^\ell 
S[-d_i]$. In this article $\exp(\A)$ indicates a \textbf{multiset}. If $\A \neq \emptyset$, then the Euler derivation $\theta_E \in D(\A)$ forms a direct summand of $D(\A)$. So we may assume that $d_1=1=\deg \theta_E$ when $\A \neq \emptyset$ is free.

Free arrangements have been a central topic in the research of hyperplane arrangements. Among them, the most important problem is so called the \textbf{Terao's conjecture} asking whether the freeness of $\A$ depends only on the intersection lattice 
$$
L(\A):=\{\cap_{H \in \B} H \mid 
\B \subset \A\}.
$$
In other words, the freeness is 
combinatorial. This is completely open, but it was shown in \cite{Z2} that 
the minimal free resolution of $D(\A)$ is not combinatorial. To approach 
Terao's conjecture, one of the main tools is Terao's addition-deletion theorem. 
For the purpose of this article, we exhibit it in a slightly different way compared as usual:
%
%

\begin{theorem}[Addition and restriction theorems, \cite{T1}]
Let $H \in \A$, $\A':=\A \setminus \{H\}$ and let $\A^H:=\{H \cap L \mid 
L \in \A'\}$. Assume that $\A'$ is free with 
$\exp(\A')=(1,d_2,\ldots,d_\ell)$. Then the following two conditions are equivalent:

\begin{itemize}
    \item [(1)]
    $\A$ is free.
    \item[(2)]
    $\A^H$ is free and $|\exp(\A') \cap \exp(\A^H)|=\ell-1$. 
\end{itemize}
In each case, for $d_i:=\exp(\A') \setminus \exp(\A^H)$, 
it holds that 
$$
\exp(\A)=(1,d_2,\ldots,d_{i-1},d_i+1,d_{i+1},\ldots,d_\ell),\ 
\exp(\A^H)=(1,d_2,\ldots,d_{i-1},d_{i+1},\ldots,d_\ell).
$$

\label{addition}
\end{theorem}

In general the condition $|\exp(\A') \cap \exp(\A^H)|=\ell-1$ above is 
described as $\exp(\A') \supset \exp(\A^H)$. However, note that $|\exp(\A') \cap \exp(\A^H)|<\ell-1$ often occurs even 
when both $\A'$ and $\A^H$ are free. So we have a first question in this article:

\begin{problem}
Assume that $\A'$ is free. Then which $\A$ makes $\A^H$ free?
\label{problem1}
\end{problem}

Also, recent developments show the following projective dimension version 
of the addition theorem. Note that $\A$ is free 
if and only if $\pd \A=0$. 

\begin{theorem}[Theorem 1.11, \cite{A9}]

(1)\,\, 
Assume that $\pd \A'=\pd \A^H=0$. Then $\pd \A \le 1$. 

(2)\,\,
Assume that $\pd \A'=0$ and $\pd \A \le 1$. Then $\pd \A^H = 0$. 
\label{pdaddition}
\end{theorem}

Now we have the second question in this article which is related to Problem \ref{problem1}:

\begin{problem}
Can we describe an algebraic structure of $D(\A)$ when $\A'$ and $\A^H$ are both free, but 
$|\exp(\A') \cap \exp(\A^H)|<\ell-1$?
\label{problem2}
\end{problem}

Explicitly, we want to know the minimal free resolution of 
$D(\A)$ under the above condition. Contrary to these problems, when $\A$ is free, 
we can describe $D(\A')$, which was proved in \cite{A5}. To see this result, let us recall the definition of 
SPOG.

\begin{define}[\cite{A5}]
We say that $\A$ is 
\textbf{SPOG} with $\POexp(\A)=(1,d_2,\ldots,d_\ell)$ and \textbf{level} $d$ if there is a 
minimal free resolution of the following:
$$
0 \rightarrow S[-d-1] \rightarrow \oplus_{i=1}^\ell S[-d_i] \oplus S[-d] 
\rightarrow D(\A) \rightarrow 0.
$$
Here $d_1=1$. For the generator $\theta_E,\theta_2,\ldots,\theta_\ell,\theta$ with 
$\deg \theta_i=d_i$ and $\deg \theta=d$ for the SPOG module $D(\A)$, $\theta_E,\theta_2,\ldots,\theta_\ell$ is called an \textbf{SPOG generator} and $\theta$ a 
\textbf{level element}. 
\label{SPOGdef}
\end{define}

It was proved in \cite{A5} (see Theorem \ref{SPOG}) that 
$\A'$ is SPOG if $\A$ is free and $\A'$ is not free. Interestingly, in this case the structure of $D(\A')$ is independent of that of $D(\A^H)$. 
However, in general $\A$ is neither free nor SPOG even if $\A'$ is free. When $\pd \A'=\pd \A^H=0$, then $\pd \A \le 1$ by Theorem \ref{pdaddition}. Also Theorem \ref{addition} shows that one additional condition for exponents say that $\pd \A=0$. 
So a weaker condition for exponents when $\pd \A'=\pd \A^H=0$ could determine the 
minimal free resolution of $D(\A)$. Namely, we can show the following, which answers Problems \ref{problem1} and \ref{problem2} partially:

\begin{theorem}
Let $\A'$ be free with $\exp(\A')=(1,d_2,\ldots,d_\ell)_\le$. 
Let $d_j < d:=d_i+d_j+|\A^H|-|\A'|  \le d_{j+1}$ for some $i<j$. 
Then the following two conditions are equivalent:
\begin{itemize}
\item[(1)]
$\A^H$ is free with $\exp(\A^H)=(1,d_2,\ldots,\hat{d}_i,\ldots,\hat{d}_j,\ldots,d_{\ell}) \cup (d)$
\item[(2)]
$\A$ is SPOG with 
$\POexp(\A)=(1,d_2,\ldots,d_i+1,\ldots,{d}_j+1,\ldots,d_\ell)$ and level $d$. 
\end{itemize}
\label{mainanother}
\end{theorem}

Theorem \ref{mainanother} can be regarded as an extension of the addition 
and restriction theorems (Theorem \ref{addition}). Namely, Theorem \ref{addition} is the 
theorem to determine the minimal free resolution of $D(\A)$ as an SPOG-module when 
$|\exp(\A') \cap \exp(\A^H)|=\ell-2$. 
The condition $d\le d_{j+1}$ is necessary, see Example \ref{eximportant} for details.

Now go back to Terao's conjecture. 
As we have seen, SPOG arrangements can be regarded as a 
close arrangements to free arrangements. 
Thus to study Terao's conjecture by using 
an inductive approach, 
it is important to study 
combinatorial dependency of SPOG arrangements. 
For that purpose, 
let us introduce the following class of arrangements.

\begin{define}
We say that $\A$ is \textbf{stair-SPOG} if 
there is $H \in \A$ such that both $\A':= \A \setminus \{H\}$ and $\A^H$ are stair-free (see Definition \ref{SF} and Theorem \ref{SFcombin}), and $\exp(\A'),\ \exp(\A^H)$ and $|\A'|-|\A^H|$ 
satisfy the conditions in Theorem \ref{mainanother}. 
Let $\mathcal{STSPOG}_\ell$ denotes the set of stair-SPOG arrangements and let 
$$
\mathcal{STSPOG}:=\bigcup_{\ell \ge 2} \mathcal{STSPOG}_\ell.
$$
\label{SPOGcombin}
\end{define}

\begin{theorem}
$\A$ is SPOG if $\A \in \mathcal{STSPOG}$. Moreover, if there is $\A,\B$ such that $\A \in \mathcal{STSPOG}$ and $L(\A) \simeq L(\B)$, then $\B$ is SPOG too.
\label{divSPOG}
\end{theorem}

The organization of this article is as follows. In \S2 we introduce several results and definitions 
for the proof of main results in this article. 
In \S3 we prove some useful results on the cardinality of the set of minimal generators. 
In \S4 we prove main results of this article. Several 
examples are also exhibited in \S4. In \S5 we study the addition theorem from the viewpoint of the cardinality of the set of minimal generators. \S6 is devoted to investigate Yoshinaga-type 
formulation of our new addition theorem. 
\medskip

\noindent
\textbf{Acknowledgements}. 
The author is
partially supported by JSPS KAKENHI Grant Number JP21H00975.
The author is grateful to Torsten Hoge for his comments to this article.

\section{Preliminaries}

In this section let us gather several definitions and results for the proof of 
main results in this article. First recall some combinatorics of arrangements. 
For the intersection lattice $L(\A)$, we can define the 
\textbf{M\"{o}bius function} 
$\mu:L(\A) \rightarrow \Z$ by $\mu(V)=1$ and by 
$$
\mu(X):=-\sum_{X \subsetneq Y \subset V,\ Y \in L(\A)} \mu(Y)
$$
for $X \in L(\A) \setminus \{V\}$. The generating function of $\mu$ is defined by 
$$
\chi(\A;t):=\sum_{X 
\in L(A)}  \mu (X) t^{\dim X}
$$
is called the \textbf{characteristic polynomial} of $\A$, which is a combinatorial invariant. 
Now let us recall several useful results on $D(\A)$.

\begin{theorem}[Terao's polynomial $B$, \cite{T1}]
Let $\C \setminus \{H\}=\C'$. Then there is a homogeneous polynomial $B 
\in S/\alpha_H S$ of degree $|\C|-1-|\C^H|$ such that 
$$
D(\C')(\alpha_H):=\{\theta(\alpha_H) \mid \theta 
\in D(\C')\} \in (\alpha_H,B).
$$
We say such $B$ a \textbf{polynomial $B$} of $(\C,H)$. 
\label{B}
\end{theorem}

\begin{theorem}[Terao's factorizatiion, \cite{T2}]
Assume that $\A$ is free with $\exp(\A)=(d_1,d_2,\ldots,d_\ell)$, then 
$$
\chi(\A;t)=\prod_{i=1}^\ell (t-d_i).
$$
\label{TF}
\end{theorem}

For $H \in \C$ and $\C':=\C \setminus \{H\}$ we have the following 
\textbf{Euler exact sequence} 
\begin{equation}
0 \rightarrow D(\C') \stackrel{\cdot \alpha_H}{\rightarrow} D(\C) 
\stackrel{\rho}{\rightarrow} D(\C^H).
\label{EE}
\end{equation}

This is not right exact in general, but it is when $\C'$ is free as follows.

\begin{theorem}[Free surjection theorem (FST), Theorem 1.13, \cite{A9}]
Let $\C=\C' \cup \{H\}$ and assume that $C'$ is free. Then 
$\rho=\rho^H:D(\C) \rightarrow D(\C^H)$ is surjective.
\label{FST}
\end{theorem}

At the end of this section we gather results on freeness and SPOGness.

\begin{theorem}[Theorem 1.4, \cite{A5}]
Let $\A$ be free with $\exp(\A)=(d_1,\ldots,d_\ell)$, $H \in \A$ and 
assume that $\A':=
\A \setminus \{H\}$ is not free. Then $\A'$ is SPOG with $\POexp(\A')=(d_1,\ldots,d_\ell)$ and level $d:=|\A'|-|\A^H|$.
\label{SPOG}
\end{theorem}

\begin{theorem}[Division theorem, Theorem 1.1, \cite{A2}]
Assume that $\A^H$ is free and $\chi(\A^H;t) \mid \chi(\A;t)$. Then $\A$ is free. Thus 
if we can show the freeness of $\A$ by using several times of the divison theorem, then 
the freeness of $\A$ is combinatorial, and such a free arrangement is called the \textbf{divisionally free arangements}. 
\label{division}
\end{theorem}

\begin{define}[Definition 4.2, \cite{A6}]
We say that $\A$ is \textbf{stair-free} if the freeness of $\A$ can be 
proved by using the addition and division theorems.
\label{SF}
\end{define}

\begin{theorem}[Theorem 4.3, \cite{A6}]
If $\A$ is stair-free, then 
its freeness depends only on $L(\A)$.
\label{SFcombin}
\end{theorem}

\section{Cardinality of the minimal set of generators}

In this section we show some new results on the cardinality of 
a minimal set of generators for $D(\A)$, which will play a key role 
to prove our main theorem. 

\begin{define}
For an arrangement $\A$, let $g(\A)$ denote 
the cardinality of 
a minimal set of generators for $D(\A)$. Clearly it is independent of the choice of the set of minimal 
generators. 
\label{ga}
\end{define} 

Moreover we define the integer for a free arrangement that measures how far 
it is from the addition.

\begin{define}
Let $\A=\A' \cup \{H\}$ and assume that $\A'$ is free. Let 
$FB(\A')$ be the set of all the basis for $D(\A')$ and for each 
$B:=\{\theta_1,\ldots,\theta_\ell\} \in FB(\A')$ define 
$$
NT(B):=|\{i\mid 1\le i \le \ell,\ 
\theta_i \not \in D(\A)\}|,
$$
and define 
$$
SNT(\A'):=\min\{NT(B) \mid B \in FB(\A')\}
$$
\label{MNT}
\end{define}

First we know the following easy facts.

\begin{lemma} 
Assume that $\A'$ is free. If $\A^H$ is free, then 
$g(\A) \le 2\ell-2$. 
\label{genless}
\end{lemma}

\noindent
\textbf{Proof}.
By Theorem \ref{FST}, 
we can choose $\theta_E,\theta_2,\ldots,\theta_{\ell-1} \in D(\A)$ as 
preimages of the basis for $D(\A^H)$. 
Let $\theta_E,\varphi_2,\ldots,\varphi_\ell$ be 
a basis for $D(\A')$. Then the Euler exact sequence (\ref{EE}) shows that 
$$
\theta_E,\theta_2,\ldots,\theta_{\ell-1},\alpha_H\varphi_2,\ldots,\alpha_H \varphi_\ell$$
generate $D(\A)$, hence $g(\A) \le \ell-1+\ell-1
=2\ell-2$. \owari
\medskip

On $g(\A)$ the following proposition is fundamental.

\begin{prop}
Let $\A'$ be free with $SNT(\A')=s$. Let 
$\theta_E,\theta_2,\ldots,\theta_\ell$ form a basis for $D(\A')$ such that 
$\theta_i \not \in D(\A)\ (2 \le i \le s+1)$ and $\theta_i  
\in D(\A)\ (i \ge s+2)$. Then $g(\A) \ge \ell+s-1$.
\label{genlowbdd}
\end{prop}

\noindent
\textbf{Proof}. 
Let $\alpha_H=x_1$, and Terao's polynomial $B \in 
S':=\K[x_2,\ldots,x_\ell]$. 
Let us show that we may choose derivations 
$\varphi_j\in D(\A)\ (j=1,\ldots,k) $ of the form $\sum_{i=2}^{s+1} f_i^j \theta_i$ such that $f_i^j \in S'$ are of positive degrees, and 
$$
\theta_E,\alpha_H \theta_2,\ldots,\alpha_H \theta_{s+1},\theta_{s+2},\ldots,\theta_\ell$$
together with derivations $\varphi_j\ (j=1,\ldots,k)$ form a minimal set of 
generators for $D(\A)$. 
Since $\deg f_i^j >0$, we know that $\theta_{s+2},\ldots,\theta_\ell,\theta_E$ cannot be removable. 
Since $f_i^j \in S'$, it is also clear that $\alpha_H \theta_2,\ldots,
\alpha_H \theta_{s+1}$ are not removable 
either. 
Assume that $\varphi_i$ is removable. We may assume that $i=1$ without loss of generality. 
Since $f_i^j \in S'$, it can be expressed as 
$$
\varphi_1=\sum_{i=2}^k g_i \varphi_i
$$
for $g_i \in S'$. Continue this among $\varphi_1,\ldots,\varphi_k$ to make them 
minimal. Then we obtain the desired minimal set of generators of the above form. 
Since their images by $\rho$ generate the rank $(\ell-1)$-module $D(\A^H)$ due to 
Theorem \ref{FST}, 
$$
|\{\theta_E,\theta_{s+2},\ldots,\theta_\ell,\varphi_1,\ldots,\varphi_k\}|=k+1
+\ell-s-1=\ell+k-s\ge \ell-1.
$$ 
So $g(\A) = \ell+k\ge \ell+s-1$.\owari
\medskip

For the rest of this section recall the following:

\begin{prop}
Let $i<j$, 
$\A'$ be free with $SNT(\A')=2$. Let $\theta_E,\theta_2,\ldots,\theta_\ell$ be a basis for 
$D(\A')$ such that $\theta_k \in D(\A)\ (k \neq i,j)$. 
Let $\theta_i(\alpha_H) =f_i\alpha_H+g_i B,\ 
\theta_j(\alpha_H)=f_j \theta_H+g_j B$ by 
Theorem \ref{B}.  Then $(g_i,g_j)=1$ and 
$D(\A)$ is generated by 
$\{\theta_k\}_{k \neq i,j} \cup\{\alpha_H \theta_i,\alpha_H \theta_j,g_j\theta_i-g_i\theta_j\}$. In particular, 
$\A$ is SPOG. 
\label{fund2}
\end{prop}

\noindent
\textbf{Proof}.
The same proof as Theorem 1.9 in  \cite{A5} shows 
that 
$\A$ is SPOG (see Proposition 3.6, \cite{Sak} too). \owari
\medskip

An immediate corollary of Proposition \ref{genlowbdd} is as follows:

\begin{cor}
Let $\A'$ be free and assume that $g(\A) =\ell+1$. Then $\A$ is SPOG.
\label{ell+1SPOG}
\end{cor}

\noindent
\textbf{Proof}.
By Proposition \ref{genlowbdd}, $g(\A) =\ell+1$ only when 
$s\le 2$. $s=0$ cannot occur, and $s=1$ implies that $\A$ is free by the addition theorem, thus 
$g(\A)=\ell$. So $s=2$ is necessary. In this case, Propotision \ref{fund2} 
shows that $\A$ is SPOG.\owari
\medskip

%

Since the freeness of $\A^H$ and $\A'$ compels that $\pd \A \le 1$ by Theorem \ref{pdaddition},  
it is natural to study which condition on $g(\A)$ makes the arrangement $\A^H$ free. 

%



\begin{theorem}
Let $\ell \ge 3$ and $\A'$ be free. 
Then $\A^H$ is free if $g(\A) \le \ell+2$.
\label{freegeneratoringeneral}
\end{theorem}

\noindent
\textbf{Proof}. By the addition and restriction theorems, Theorem \ref{FST}, 
Corollary \ref{ell+1SPOG} and the explicit form of SPOG generators as in Proposition 
\ref{fund2}, the statement follows if 
$g(\A) \le \ell+1$. Assume that $g(\A)=\ell+2$. 
Then by Proposition \ref{genlowbdd}, we have a basis $\theta_E,\theta_2,\ldots,\theta_\ell$ for $D(\A')$  such that 
$\theta_i \not \in D(\A)$ for $i 
\ge \ell-2$ and $\theta \in D(\A)$ for $i \le \ell-3$. We may 
assume that $\alpha_H=x_1$. Also for the polynomial $B$ of $(\A,H)$, we may assume that $B \in \K[x_2,\ldots,x_\ell]$ and 
$\theta_i(x_1)=f_i B$ modulo $x_1$ and $f_i \in \K[x_2,\ldots,x_\ell]$ as in the proof of 
Proposition \ref{genlowbdd}. 
Then clearly $D(\A)$ has a minimal set of generators with cardilality $\ell+2$ of the form 
$$
\theta_E,\theta_2,\ldots,\theta_{\ell-3},x_1\theta_{\ell-2},
x_1\theta_{\ell-1},x_1\theta_{\ell},\varphi_1,\varphi_2,
$$
where $\varphi_j$ is of the linear combination of $\theta_{\ell-2},\theta_{\ell-1},\theta_\ell$ over $\K[x_2,\ldots,x_\ell]$. By 
Theorem \ref{FST}, the images of $\theta_E,
\theta_2,\ldots,\theta_{\ell-3},\varphi_2,\varphi_3$ by $\rho$ have to generate $D(\A^H)$. Since $\mbox{rank}_{S/\alpha_H S} D(\A^H)=\ell-1$, 
it holds that $\A^H$ is free. \owari
\medskip

By using results above, we can show the following 
proposition which is fundamental on the relation between free and SPOG arrangements.

\begin{prop}
Let $\A$ be SPOG, $H \in \A$ and $\A':=
\A \setminus \{H\}$. 
If $\A'$ is free, then there is 
SPOG generators $\theta_1=\theta_E,\theta_2,\ldots,\theta_\ell$, a level element  $\varphi$ and two distinct integers 
$1 < s < t \le \ell$ such that 
$$
\theta_E,\theta_2,\ldots,\theta_{s-1},
\theta_s/\alpha_H,\theta_{s+1},\ldots,
\theta_{t-1},
\theta_t/\alpha_H,\theta_{t+1},\ldots,\theta_\ell
$$
form a free basis for $D(\A')$.
\label{fund}
\end{prop}

\noindent
\textbf{Proof}.
By Proposition \ref{genlowbdd} and the 
assumption that $g(\A)=\ell+1$, for the basis 
$\theta_E,\varphi_2,\ldots,\varphi_\ell$ for $D(\A')$, we may assume that 
$\varphi_i \in D(\A)\ (i \le \ell-2)$ and $\varphi_{\ell-1},\varphi_\ell 
\not \in D(\A)$. Then 
$$
\theta_E,\varphi_2,\ldots,\varphi_{\ell-2},\alpha_H \varphi_{\ell-1}, 
\alpha_H\varphi_\ell,f_{\ell}\varphi_{\ell-1}-f_{\ell-1}
\varphi_\ell
$$
form a minimal set of generators for $D(\A)$, where 
$\varphi_i(\alpha_H)=f_i B$ modulo $\alpha_H$ for $i=\ell-1,\ell$ and 
$B$ is Terao's polynomial, which is the required generators for $D(\A)$.\owari
\medskip

\section{Proof of main results}


In this section let us prove Theorems \ref{mainanother} and \ref{divSPOG}.
\medskip

\noindent
\textbf{Proof of Theorem \ref{mainanother}}. 
First prove $(2) \Rightarrow (1)$. 
It is easy to see that, for the basis $\theta_E,\varphi_2,\ldots,\varphi_\ell$ for 
$D(\A')$ with 
$\deg \varphi_k=d_k$, the SPOG generator for $D(\A)$ is of the form 
$\{\theta_k\}_{k \neq i,j} \cup \{\alpha_H \theta_k\}_{k=i,j} \cup \{f_j\theta_i-f_i\theta_j\}$ by 
Proposition \ref{fund2}. 
Thus Theorem \ref{FST} shows that $D(\A^H)$ is generated by the image of 
$\{\theta_k\}_{k \neq i,j} \cup\{f_j\theta_i-f_i\theta_j\}$. Since $\mbox{rank}_S D(\A^H)=\ell-1$, it follows that $\A^H$ is free with 
given 
exponents.

Conversely, prove $(1) \Rightarrow (2)$. 
Assume that 
$\A^H$ is free with the given exponents above. 
In this assumption, Terao's addition theorem shows that $\A$ is not free since 
$\exp(\A^H) \not 
\subset \exp(\A')$. 
By Theorem \ref{FST}, the Euler restriction map 
$\rho^H:D(\A) \rightarrow D(\A^H)$ is surjective. 
Thus there are $S$-independent 
derivations $\theta_E,\theta_2,\ldots,\hat{\theta}_i,\ldots,\hat{\theta}_j,\ldots,\theta_\ell,\theta 
\in D(\A)$ such that $\deg \theta_k=d_k,\ \deg \theta=d$, and their images by $\rho^H$ form a basis for 
$D(\A^H)$. Let $\theta_E,\varphi_2,\ldots,\varphi_\ell$ be a basis for $D(\A')$ with 
$\deg \varphi_k=d_k$. We may assume that $d_i<d_{i+1}$. By the usual 
degree argument, we may assume that $\varphi_k=\theta_k$ for 
$k<i$. Next, assume that $\varphi_s=\theta_s$ for all $i\neq s<k$, where 
$k>i$ and $d_i< d_k < d_j$. 
Since $\theta_k \in D(\A) \subset D(\A')$, we have an expression 
$$
\theta_k=\sum_{s \neq i,s=1}^{k-1} f_s \theta_s+\sum_{s=k}^\ell f_s \varphi_s+f_i \varphi_i.
$$
Note that $f_t=0$ if $d_t>d_k$. Assume that there is $t$ with $d_t=d_k$ and $f_t \neq 0$. 
Then first, take all the sum of $0 \neq \sum f_s \varphi_s$ for all $s$ with $d_s=d_k,\ f_s\neq 0$ and 
let rename this derivation as $\varphi_k$. Then replacing $\varphi_k$ by $\varphi_k+\sum_{s \neq i,s=1}^{k-1} f_s \theta_s +f_i \varphi_i =\theta_k \in D(\A')$, 
we may assume that 
$\theta_k=\varphi_k$. So let us assume that $f_t=0$ for all $t$ with 
$d_t=d_k$. 
If $f_i=0$, then $\theta_k$ is a linear combination of $\theta_E,\theta_2,\ldots,\theta_{k-1}$, a contradiction. 
So $f_i \neq 0$ by the independency. 
Assume that $\varphi_i \in D(\A)$. 
Then $\rho^H(\varphi_i)$ can be expressed as a linear combination of 
$\rho^H(\theta_E),\rho^H(\theta_2),\ldots,\rho^H(\theta_{i-1})$ by 
the 
reason of degrees. Since $D(\A^H)$ is free, we may assume that $\alpha_H \mid \varphi_i$, contradicting the minimality of the generator of $\theta_E,\varphi_2,\ldots,\varphi_\ell$. 
So $\varphi_i \not \in D(\A)$. 
Since $\theta_s=\varphi_s \in D(\A)$ if $i \neq s<k$,
it holds that $f_i \varphi_i \in D(\A)$. So $\alpha_H \mid f_i$.
Thus the image of $\theta_k$ by $\rho^H$ is again a linear combination of 
the images of $\theta_E,\theta_2,\ldots,\theta_{k-1}$, a contradiction. 
So we may assume that $\theta_k=\varphi_k \in D(\A)$ for all $k$ with $d_i<d_k<d_j$. 

We may assume that $d_j<d_{j+1}$ or $d=\ell$. Now replace $\theta_k$ by $\varphi_k$ for $k<j$ satisfying $d_k=d_j$ as far
as possible. 
Assume that 
$$
\theta_k=\sum_{s \neq i,s=1}^{k-1} f_s \theta_s+\sum_{s=k}^\ell f_s \varphi_s+f_i \varphi_i
$$
and $f_s=0$ for all $s$ with $d_s=d_j=d_k$. Then the same argument as in the previous paragraph 
implies that $\alpha_H \mid f_i$, and thus $\rho^H(\theta_k)$ 
cannot be a part of basis for $D(\A^H)$. 
Hence 
first, take all the sum of $0 \neq \sum f_s \varphi_s$ for all $s$ with $d_s=d_k=d_j,\ f_s\neq 0$ and 
let rename this non-zero 
derivation as $\varphi_k$. Then replacing $\varphi_k$ by $\varphi_k+\sum_{s \neq i,s=1}^{k-1} f_s \theta_s +f_i \varphi_i =\theta_k \in D(\A)$, 
it holds that $\varphi_k=\theta_k$ for $j+1>k \neq i,j$. 
Note that $\varphi_j \not \in D(\A)$ 
by the same 
proof of $\varphi_i \not \in D(\A)$. 

Next express $\theta$ in the following form by using the fact that 
$d =\deg \theta <d_{j+1}$:
\begin{equation}
\theta=\sum_{i\neq k=1}^{j-1} f_k \theta_k+f_i \varphi_i+f_j \varphi_j.
\label{eq1}
\end{equation}
By the independency of images of $\theta_k$ by $\rho$, at least one of $f_i,f_j$ is not zero. Assume that only one of 
them is not zero, say $f_i\neq 0$ and $f_j=0$. 
Then (\ref{eq1}) combined with the fact that $\varphi_i \not \in D(\A)$ 
shows that $\alpha_H \mid f_i$. 
By replacing $\theta$ by $\theta-\sum_{i\neq k=1}^{j-1} f_k$, we may assume that 
$\theta=f_i \varphi_i$ with $\alpha_H \mid f_i$. So $\rho^H(\theta)=0$, a contradiction. 
%
Thus both $f_i$ and $f_j$ are not zero.
So replacing $\theta$ by $\theta-\sum_{i\neq k=1}^{j-1} f_k \theta_k$, we may assume that 
$$
\theta=f_i\varphi_i+f_j \varphi_j.
$$
Recall 
that 
$\varphi_i,\varphi_j \not \in D(\A)$. By using Theorem \ref{B} and Proposition \ref{fund2}, letting $\varphi_i(\alpha_H)=g_i B$ 
and $\varphi_j(\alpha_H)=g_j B$ modulo $\alpha_H$, 
$(g_i,g_j)=h$ and $g_i=hh_i,\ g_j=-hh_j$, it holds that 
$\varphi:=h_j\varphi_i-h_i\varphi_j \in D(\A)$. 
By the construction, $\deg \varphi\le d=\deg \theta=d_i+d_j-|\A'|+|\A^H|=d_i+d_j-\deg B$. 
Assume that $\deg 
\varphi < d$. 
Then $h_i\varphi_j-h_j\varphi_i \in D(\A)$, so its image by $\rho$ is expressed as a linear combination of 
$\{\rho(\theta_k)\}_{i.j \neq k=1}^{j-1}$. Thus 
$$
h_i\varphi_j-h_j\varphi_i-\sum_{i,j \neq s=1}^{j-1} a_s \theta_s\in 
\alpha_H D(\A')_{<d_{j+1}}=\alpha_H\langle 
\{\theta_k\}_{i,j \neq k}^{j-1} \cup\{\varphi_i,\varphi_j\} \rangle_S
$$
for some $a_s \in S$. 
So $\alpha_H \mid h_i$ and  $\alpha_H \mid h_j$, 
contradicting $\varphi_i 
\not \in D(\A)$. Thus we may assume that $h=1$ and 
$\deg (g_i\varphi_j-g_j \varphi_i)=d=\deg \theta$.
Now express 
$$
\theta=\sum_{i,j \neq k=1 }^{\ell} a_k \theta_k+a_i \varphi_i+
a_j \varphi_j.
$$
First assume that $d < d_{j+1}$. Then the above expression is 
$$
D(\A) \ni \theta=\sum_{i,j \neq k=1 }^{j-1} a_k \theta_k+a_i \varphi_i+
a_j \varphi_j.
$$
Then 
Proposition \ref{fund2} shows that $a_i\varphi_i+a_j\varphi_j$ is of the form 
$x_1(b_i\varphi_i+b_j\varphi_j)+c(g_j\theta_i-g_i\theta_j)$, where $\alpha_H=x_1$. 
If $c=0$, then the same argument as the above shows that $\rho(\theta)=0$, a contradiction. Thus 
$c\neq 0$, and
$\deg \theta=\deg (g_j\theta_i-g_i \theta_j)$ shows that $c=1$. So we may assume that 
$\theta=g_i\theta_j-g_j\theta_i$.

Next assume that $d=d_{j+1}$. Then the above expression is 
$$
D(\A) \ni \theta=\sum_{i,j \neq k=1 }^{j-1} a_k \theta_k+\sum_{k=j+1}^p a_k \varphi_k+a_i \varphi_i+
a_j \varphi_j.
$$
Here $d=d_{j+1}=\cdots=d_p<d_{p+1}$ and $a_k \in \K$ for $j+1\le k \le p$. If some $a_k\neq 0$ for 
$j+1\le k \le p$, then 
we may assume, by the same argument as 
the above, that $\theta=\varphi_k$ for some  $j+1\le k \le p$. Clearly this replacement continues 
in the range of $j+1\le k \le p$ in such a way that $\theta_k=
\varphi_{k+1}$ for 
$j+1\le k \le p-1$. Thus in the final step, we have 
$$
D(\A) \ni \theta_p=\sum_{i,j \neq k=1 }^{p-1} a_k \theta_k+a\theta+a_i \varphi_i+
a_j \varphi_j.
$$
Then the same argument shows that both of $a_i$ and $a_j$ are not zero, and the independency, 
$d_{p}=d<d_{p+1}$ and 
Proposition \ref{fund2} shows that $a_i \varphi_i+
a_j \varphi_j=x_1(b_i\varphi_i+b_j\varphi_j)+c(g_j\theta_i-g_i\theta_j)$. 
Now exchange $\theta_p$ and $\theta$ 
to obtain that $\theta=g_j\theta_i-g_i\theta_j$. 

Now continue the replacement of $\theta_k$ by $\varphi_k$ for $k>j$. Assume that 
$\theta_k=\varphi_k$ for $ i,j 
\neq k<s$. Let 
\begin{equation}
\theta_s=\sum_{i,j\neq k=1}^{\ell} a_k \theta_k+a_i \alpha_H \varphi_i+a_j \alpha_H \varphi_j+a\theta.
\label{eq2}
\end{equation}
If $a_k=0$ for all $k$ 
with $d_k=d_s$, then 
replacing 
$\theta_s$ by 
$\theta_s-\sum_{i,j\neq k=1}^{\ell} a_k \theta_k-a_i \alpha_H \varphi_i-a_j \alpha_H \varphi_j,
$
we may assume that $\theta_s=a\theta$, contradicting the independency of 
$\rho(\theta_E),\{\rho(\theta_k)\}_{i,j \neq k=2}^\ell, \rho(\theta)$. 
So by the same replacement we may assume that $\theta_s=\varphi_s$ up to $s = \ell$. Now it holds that $\varphi_k=\theta_k \in D(\A)$ unless 
$k=i,j$. So Proposition \ref{fund2} shows that 
$\A$ is SPOG with $\POexp(\A)=(1,d_2,\ldots,d_i+1,\ldots,d_j+1,\ldots,d_\ell)$ and level $d=d_i+d_j-|\A'|+|\A^H|$. \owari
\medskip



The following case is the most practical to apply Theorem \ref{mainanother}.

\begin{cor}
Let $\A'$ be free with $\exp(\A')=(1,d_2,\ldots,d_\ell)_\le$ and $d:
=d_i+d_\ell+|\A^H|-|\A'| > d_\ell$ for some $i$. 
Then the following two conditions are equivalent:
\begin{itemize}
    \item [(1)]
$\A^H$ is free with $\exp(\A^H)=(1,d_2,\ldots,\hat{d}_i,\ldots,d_{\ell-1}) \cup (d)$
\item[(2)]
$\A$ is SPOG with 
$\POexp(\A)=(1,d_2,\ldots,d_{i-1},d_i+1,d_{i+1},\ldots,d_{\ell-1},d_\ell+1)$ and level $d$. 
\end{itemize}
\label{onedirection}
\end{cor}

\noindent
\textbf{Proof}.
Clear by the proof of Theorem \ref{mainanother}.
\owari
\medskip

Let us apply Theorem \ref{mainanother} to some examples.

\begin{example}
Let $\A$ be the Weyl arrangement of the type $A_4$ defined by 
$$
Q(\A)=\prod_{i=1}^4 x_i \prod_{1 \le i < j \le 4} (x_i-x_j)=0.
$$
$\A$ is well-known to be free with $\exp(\A)=(1,2,3,4)$. Let $\A \not \ni 
H:x_1-x_2+2x_3-2x_4=0$ and let $
\B:= \A 
\cup \{H\}$. Then $|\B^H|=9 <10=|\A|$. 
It is easy to show that 
$\B^H$ is free with $\exp(\B^H)=(1,4,4)$. Since 
$$
2+3-|\A|+|\B^H|=5-1=4 \ge 4,
$$  
Theorem \ref{mainanother} shows that 
$\B$ is SPOG with $\POexp(\B)=(1,3,4,4)$ and level $4$.

Note that $Q(\B^H)$ is
$$
x_2x_3x_4(x_2-x_3)(x_2-x_4)(x_3-x_4)(x_2-2x_3+2x_4)(x_2-3x_3+2x_4)(x_2-2x_3+x_4).
$$
Let $L:x_2=0$ and ket $\C:=\B^{H\cap L}$. Then it is easy to show that 
$\chi(\C;t)=(t-1)(t-4)$ and $\chi(\B^H;t)=(t-1)(t-4)^2$. Thus $\B^H$ is divisionally free as in 
Theorem \ref{division}. Since 
$\A$ is divisionally free too, by Theorem \ref{TF}, the freeness and 
exponents of $\A$ and $\B^H$ are both combinatorial. 
Thus Theorem \ref{mainanother} shows that the SPOGness of 
$\B$ is combinatorially determined.
\label{ex1}
 \end{example}


\noindent
\textbf{Proof of Theorem \ref{divSPOG}}.
Clear by Theorem \ref{mainanother} and Theorem \ref{SFcombin}.\owari
\medskip

We can use Theorem \ref{mainanother} to show the combinatorial freeness of 
arrangements by using a non-free but SPOG arrangements. Let us check it by the following example:

\begin{example}
Let $\A$ be the Weyl arrangement of the type $B_4$ defined by 
$$
Q(\A)=\prod_{i=1}^4 x_i \prod_{1 \le i<j \le 4} (x_i^2-x_j^2)=0.
$$
Let  
$\A':=\A \setminus \{H\}$, where $H:x_1=0$. Let $L:x_1+x_2+x_3=0$. 
We know that $\A'$ is divisionally free with $\exp(\A')=(1,3,5,6)$, and 
$(\A'\cup\{L\})^L$ is also divisionally free with $\exp(\A^H)=(1,5,7)$. Thus Theorem \ref{mainanother} 
confirms that $\B:=\A'\cup\{L\}$ is SPOG with $\POexp(\B)=(1,4,5,7)$ and level $7$. 
Next let $\C:=\B \cup \{H\}$. Then $|\C^H|=9$, so $|\B|-|\C^H|=16-9=7$. 
Thus the generators of degree $1,4,5$ for $D(\B)$ are in $D(\C)$ too, and we may assume that one of 
two generators of degree $7$ is in $D(\C)$ by Theorem \ref{B}. Since the relation in $D(\B)$ are among three derivations of degrees $4,7,7$ by the explicit construction of 
SPOG generators and a level element in Proposition \ref{fund2}, we know that these 
$4$-basis elements in $D(\C)$ are $S$-independent. So $\C$ is combinatorially free with $\exp(\C)=(1,4,5,7)$ since the SPOGness of $\B$ and $|\C|$ are both combinatorial.
\label{deletedB}
\end{example}

Unfortunately, there are cases in which $\A'$ and $\A^H$ are free, $\A$ SPOG but Theorem \ref{mainanother} cannot be applied.

\begin{example}

Let $\A$ be the Weyl arrangement of the type $\A_4$ and let 
$H:x_1+x_2+x_3=0$. Say that $\B:=\A \cup \{H\}$. Then $\B^H$ is free with 
$\exp(\B^H)=(1,4,5) \not \subset (1,2,3,4)=\exp(\A)$. We can check that 
$\B$ is SPOG with $\POexp(\B)=(1,3,4,4)$ and level $5$ by using Macaulay2 in \cite{GS}, but we cannot apply Theorem \ref{mainanother}.

\label{ex2}
 \end{example}

\begin{problem}
Generalize Theorem \ref{mainanother} to all the case when $\A'$ and $\A^H$ are free, 
$\ell-2 \le |\exp(\A') \cap \exp(\A^H)| \le \ell-1$ and $\A$ not free. More explicitly, could it be possible that 
$\A',\ \A^H$ are both free but $\A$ not SPOG? Even in this case $\pd_S D(\A) =1$. 
\label{mainanotherdiff}
\end{problem}

In fact, the condition $d \le d_{j+1}$ is necessary. Let us see the following example.

\begin{example}
Let $\A'$ be an arrangement in $\R^4$ 
defined by 
$$
Q(\A')=(x_1+x_2+x_3+x_4)\prod_{i=1}^4 x_i \prod_{i=2}^4 (x_1+x_i) \prod_{i=2}^4
(x_1+x_2+x_3+x_4-x_i).
$$
Then $\A'$ is free with $\exp(\A')=(1,3,3,4)$. 
Let $H=x_2+x_3+x_4$ and let $\A:=
\A' \cup\{H\}$. Then $\A^H$ is free with $\exp(\A^H)=(1,4,5)$. So 
$\exp(\A') \cap \exp(\A^H)=\{1,4\}$. However, Macaulay2 in 
\cite{GS} shows that $D_0(\A)$ is not SPOG but has a following minimal 
free resolution:
$$
0 \rightarrow 
S[-5]\oplus S[-6] 
\rightarrow 
S[-4]^3 \oplus S[-5]^2 \rightarrow D_0(\A) 
\rightarrow 0.
$$
Since $d=5=3+3-|\A'|+|\A^H|>4$, 
the condition $d\le d_{j+1}$ in Theorem \ref{mainanother} is necessary.
\label{eximportant}
\end{example}

If we remove the assumption on $\exp(\A')$ and 
$\exp(\A^H)$ in Theorem \ref{mainanother}, 
we have an example related to Problem \ref{mainanotherdiff}.

\begin{example}
Let 
$$
Q(\A')=\prod_{i=1}^4 x_i \prod_{i=1}^3 (x_i^2-x_4^2)(x_i^2-4x_4^2)
\prod_{i=2}^3 (x_i^2-9x_4^2)
(x_3^2-16x_4^2).
$$
Then $\A'$ is free with $\exp(\A')=(1,5,7,9)$. Let 
$H_1:x_2+x_3+7x_4=0$ and $H_2:x_1+x_2+x_3=0$ and let 
$\A_i:=\A' \cup\{H_i\}$. Then $\A_1$ is SPOG with $\POexp(\A_1)=(1,5,8,10)$ and level $15$. So Theorem \ref{mainanother} shows that 
$\A^H$ is free with exponents $(1,5,15)$ and vice versa. On the other hand, $\A_2^{H_2}$ is free with exponents 
$(1,10,11)$, and $D(\A_2)$ has a minimal free resolution 
$$
0 \rightarrow S[-11]\oplus S[-12] 
\rightarrow 
 S[-6]\oplus S[-8] \oplus S[-10]^2 \oplus S[-11] 
 \rightarrow D_0(\A_2) \rightarrow 0.
 $$
 So in general exponents case, it could occur that $\A$ is of projective dimension one, not SPOG, $\A'$ and $\A^H$ are free. 
 
Note that the freeness of $\A^H$ follows from $g(\A) \le 6$ and the 
freeness of $\A'$ by Theorem \ref{freegeneratoringeneral}.
 \end{example}

\end{document}